\documentclass{article}
\usepackage{amsmath,amssymb}
\newtheorem{df}{Definition}[section]
\newtheorem{thm}[df]{Theorem}
\newtheorem{lem}[df]{Lemma}
\title{Uncertainty of Poisson wavelets}
\author{Ilona Iglewska-Nowak\footnote{West Pomeranian University of Technology in Szczecin, School of Mathematics, al. Pias\-t\'ow 17, 70--310 Szczecin, Poland}}

\begin{document}

\maketitle

\bibliographystyle{amsplain}

\begin{abstract}Poisson wavelets are a powerful tool in analysis of spherical signals. In order to have a deeper characterization of them, we compute their uncertainty product, a quantity introduced for the first time by Narcowich and Ward in~\cite{NW96} and used to measure the trade-off between the space and frequency localization of a function. Surprisingly, the uncertainty product of Poisson wavelets tends to the minimal value in some limiting cases. This shows that in the case of spherical functions, not only Gauss kernel has this property.

\end{abstract}

\begin{bfseries}Key words and phrases:\end{bfseries}  uncertainty product, Poisson multipole wavelets, spherical wavelets \\
\begin{bfseries}2010 Mathematics Subject Classification Number:\end{bfseries} 42C40

\section{Introduction}

The uncertainty principle is a very important result in quantum mechanics. It can be also studied independently of physical interpretations, as a property of square integrable functions. In 1990s an uncertainty principle has been proven for functions over the sphere~\cite{NW96,RV97}. The so-called uncertainty product can be understood as a measure for sharpness of a function, simultaneously in space and frequency. In the case of analyzing functions such as wavelets it yields an information of the resolution of the transform that can be obtained. For this reason, we compute the uncertainty product of spherical wavelets most frequently used in practise, namely Poisson wavelets.

Poisson wavelets, introduced in~\cite{HCM03}, have proven to be very useful for their numerous practical properties~\cite{CPMDHJ05,HI07}. There exist explicit expressions of Poisson wavelets in terms of spherical variables~\cite{IH10} whereas most other wavelet families are given as series of (hyper)spherical harmonics~\cite{FGS-book,FW-C}, the wavelets are well localized in space, and there exist discrete frames of them~\cite{IH10,IIN16WF}. Moreover, they have Euclidean limit property, i.e., for small scales they behave like wavelets over the Euclidean plane~$\mathbb R^n$~\cite{IIN14CWT,IIN15PW}.

The paper is organized as follows. In Section~\ref{sec:sphere} we introduce Poisson wavelets as well as notions and concerned with the spherical uncertainty principle. Section~\ref{sec:uncertainty} is devoted to the computation of the uncertainty product of Poisson wavelet and discussion of the result. Some elementary but tedious computations are postponed to Appendix.

\section{Preliminaries}\label{sec:sphere}

Denote by~$\zeta$ point $(r,0,\dots,0)$, $r=e^{-\rho}$, $\rho\in(0,\infty)$, inside the unit sphere $\mathcal S^n=\{x\in\mathbb R^{n+1}:\,|x|=1\}$, and let $x\in\mathcal S^n$. Poisson kernel is given by
$$
p_\rho(x)=\frac{1}{\Sigma_n}\frac{1-|\zeta|^2}{|\zeta-x|^{n+1}}=\frac{1}{\Sigma_n}\frac{1-\rho^2}{(1-2\rho\cos\vartheta+\rho^2)^{(n+1)/2}},
$$
where
$$
\Sigma_n=\int_{\mathcal S^n}d\sigma=\frac{2\pi^{\lambda+1}}{\Gamma(\lambda+1)}
$$
is the Lebesgue measure of the sphere,
$$
\lambda=\frac{n-1}{2},
$$
and
$$
r\cos\vartheta=\zeta\cdot x,
$$
i.e., $\vartheta$ is the first spherical variable of~$x$. Poisson wavelet of order~$m$, $m\in\mathbb N$, at a scale~$\rho$ is given recursively by
\begin{align*}
g_\rho^1&=\rho r\,\partial_r\,p_\rho,\\
g_\rho^{m+1}&=\rho r\,\partial_r\,g_\rho^m,
\end{align*}
see \cite[Definition~3.1]{IIN15PW}. Lemma~3.2 in \cite{IIN15PW} states that the Gegenbauer expansion of Poisson wavelets, i.e., the representation of a rotation-invariant function in terms of Gegenbauer polynomials~$\mathcal C_l^\lambda$, is given by
$$
g_\rho^{m}(x)=\frac{1}{\Sigma_n}\,\sum_{l=0}^\infty\frac{l+\lambda}{\lambda}\,(\rho l)^m e^{-\rho l}\,\mathcal C_l^\lambda(\cos\vartheta).
$$

According to~\cite{RV97}, for a twice differentiable spherical function~$f$ with $\int_{\mathcal S^n}x\,|f(x)|^2\,dx\ne0$, its variances in space and momentum domain are given by
$$
\text{var}_S(f)=\left(\frac{\int_{\mathcal S^n}|f(x)|^2\,d\sigma(x)}{\int_{\mathcal S^n}x\,|f(x)|^2\,d\sigma(x)}\right)^2-1
$$
and
$$
\text{var}_M(f)=-\frac{\int_{\mathcal S^n}\Delta^\ast f(x)\cdot \bar f(x)\,d\sigma(x)}{\int_{\mathcal S^n}|f(x)|^2\,d\sigma(x)},
$$
where $\Delta^\ast$ is Laplace--Beltrami operator on~$\mathcal S^n$, and the quantity
$$
U(f)=\sqrt{\text{var}_S(f)}\cdot\sqrt{\text{var}_M(f)}
$$
is called uncertainty product of~$f$. \cite[Theorem~1.2]{RV97} states that for rotation-invariant functions $f\in\mathcal L^2(\mathcal S^n)\cap\mathcal C^2(\mathcal S^n)$,
\begin{equation}\label{eq:UP}
U(f)\geq\frac{n}{2},
\end{equation}
and the lower bound is optimal. The inequlatiy~\eqref{eq:UP} is the uncertainty principle for spherical functions.

The variances in space and momentum domains for zonal functions can be computed from their Gegenbauer coefficients, see \cite[Lemma~4.2]{IIN16MRes}.

\begin{lem}\label{lem:varS_varM} Let a zonal $\mathcal L^2(\mathcal S^n)$--function be given by its Gegenbauer expansion
$$
f(t)=\sum_{l=0}^\infty\widehat f(l)\,C_l^\lambda(t).
$$
Its variances in space and momentum domain are equal to
\begin{align}
\text{var}_S(f)&=\left(\frac{\sum_{l=0}^\infty\frac{\lambda}{l+\lambda}\,\binom{l+2\lambda-1}{l}\,|\widehat f(l)|^2}{\sum_{l=0}^\infty\binom{l+2\lambda}{l}\,
   \frac{\lambda^2\left[\overline{\widehat f(l)}\,\widehat f(l+1)+\widehat f(l)\,\overline{\widehat f(l+1)}\right]}{(l+\lambda)(l+\lambda+1)}}\right)^2-1,\label{eq:varS}\\
\text{var}_M(f)&=\frac{\sum_{l=1}^\infty\frac{ l\lambda(l+2\lambda)}{l+\lambda}\,\binom{l+2\lambda-1}{l}\,|\widehat f(l)|^2}
   {\sum_{l=0}^\infty\frac{\lambda}{l+\lambda}\,\binom{l+2\lambda-1}{l}\,|\widehat f(l)|^2},\label{eq:varM}
\end{align}
whenever the series are convergent.
\end{lem}

\section{Uncertainty}\label{sec:uncertainty}

\begin{thm}
Let $g_\rho^m$ be Poisson wavelet of order~$m\geq1$ over~$\mathcal S^n$, $n\geq2$. Its variance in space domain behaves in limit $\rho\to0$ like
$$
\text{var}_S(g_\rho^m)=\begin{cases}\frac{n^2-3n+2m+2}{(n+2m-1)(n+2m-2)}\cdot\rho^2-\frac{4\,(n-1)^2(n-3)\,m}{(n+2m-1)^2(n+2m-2)(n+2m-3)}\cdot\rho^3
   &\hspace{-1em}+\mathcal{O}(\rho^4),\\
&\hspace{-3em}n\geq5,\,n=2,\\
\frac{(m+3)}{2m^2+5m+3}\cdot\rho^2-\frac{2m\,(4m^2+2m+21)}{3\,(2m+3)^2(2m^2+3m+1)}\cdot\rho^3+\mathcal O(\rho^4),&n=4,\\
\frac{1}{2m+1}\cdot\rho^2-\frac{2\,(m-1)(m+1)(m+5)}{3\,(2m+1)}\cdot\rho^3+\mathcal O(\rho^4),&n=3.\end{cases}
$$
The variance in momentum domain behaves like
$$
\text{var}_M(g_\rho^m)=\frac{(n+2m)(n+2m+1)}{4\rho^2}+\frac{(n-1)\,m\,(n+2m)}{(n+2m-1)\,\rho}+\mathcal O(1),\qquad\rho\to0,
$$
and the uncertainty product is given by
$$
U(g_\rho^m)=\begin{cases}\frac{1}{2}\,\sqrt{\frac{(n+2m)(n+2m+1)[n^2-3n+2\,(m+1)]}{(n+2m-1)(n+2m-2)}}
   -\frac{(n-1)\,m\,[3n^2-4n\,(m+3)-4m^2+8m+9]}{n+2m-3}&\\
\quad\cdot\sqrt{\frac{n+2m}{(n+2m+1)[n^2-3n+2(m+1)](n+2m-1)^3(n+2m-2)}}\cdot\rho+\mathcal O(\rho^2),&\hspace{-3em}n\geq5,\,n=2,\\
\sqrt{\frac{(m+3)(m+2)(2m+5)}{2(m+1)(2m+3)}}\cdot\left[1-\frac{m\,(8m^3-12m^2-74m+51)}{3\,(m+3)(2m+1)(2m+3)(2m+5)}\cdot\rho\right]+\mathcal O(\rho^2),&n=4,\\
\sqrt{\frac{(m+2)(2m+3)}{2(2m+1)}}\cdot\left[1-\frac{m^5+8m^4+16m^3+2m^2-20m-10}{3\,(m+1)(m+2)}\cdot\rho\right]+\mathcal O(\rho^2),&n=3,\end{cases}
$$
for $\rho\to0$.
\end{thm}

\begin{bfseries}Proof.\end{bfseries} Consider the function
$$
f_\rho^m(x)=\frac{\Sigma_n}{\rho^m}\cdot g_\rho^m(x),\quad x\in\mathcal S^n.
$$
Its Gegenbauer coefficients are given by
$$
\widehat{f_\rho^m}(l)=\frac{l+\lambda}{\lambda}\,l^me^{-\rho l},\qquad l\in\mathbb N_0.
$$
According to Lemma~\ref{lem:varS_varM}, the space and momentum variances of $g_\rho^m$ are equal to the space and momentum variances of $f_\rho^m$, respectively. We shall use the notation~$S_m$ for the expression
$$
S_m(\rho)=\sum_{l=0}^\infty\binom{l+2\lambda-1}{l}\,l^me^{-2\rho l};
$$
then, the variances in space and momentum domain of~$f_\rho^m$ can be written as
\begin{align}
\text{var}_S(f_\rho^m)
   &=\left(\frac{\frac{1}{\lambda}\,S_{2m+1}+S_{2m}}{e^{-\rho}\,\sum_{j=0}^m\binom{m}{j}\left(\frac{1}{\lambda}\,S_{m+j+1}+2S_{m+j}\right)}\right)^2-1,\label{eq:var_S}\\
\text{var}_M(f_\rho^m)&=\frac{\frac{1}{\lambda}\,S_{2m+3}+3S_{2n+2}+2\lambda\,S_{2n+1}}{\frac{1}{\lambda}\,S_{2m+1}+S_{2m}}.\notag
\end{align}
(Since Gegenbauer coefficients of~$f_\rho^m$ are real, the denominator in~\eqref{eq:var_S} is obtained from the one in~\eqref{eq:varS} in the following way:
\begin{align*}
\sum_{l=0}^\infty&\binom{l+2\lambda}{l}\frac{2\lambda^2\widehat f(l)\widehat f(l+1)}{(l+\lambda)(l+\lambda+1)}
   =\sum_{l=0}^\infty2\,\binom{l+2\lambda}{l}\,l^me^{-\rho l}\cdot(l+1)^me^{-\rho(l+1)}\\
&=\sum_{l=0}^\infty\frac{2(l+2\lambda)}{2\lambda}\binom{l+2\lambda-1}{l}\,l^me^{-(2\rho l+1)}\sum_{j=0}^m\binom{m}{j}l^j\\
&=e^{-\rho}\sum_{l=0}^\infty\binom{l+2\lambda-1}{l}\,e^{-2\rho l}\sum_{j=0}^m\binom{m}{j}\left(\frac{l}{\lambda}+2\right)l^{m+j}\\
&=e^{-\rho}\sum_{j=0}^m\binom{m}{j}\sum_{l=0}^\infty\binom{l+2\lambda-1}{l}\left(\frac{l^{m+j+1}}{\lambda}+2l^{m+j}\right)\,e^{-2\rho l}.)
\end{align*}
It can be proven by induction with respect to $n=2\lambda+1$ that
$$
S_0(\rho)=\frac{1}{(1-e^{-2\rho})^{n-1}}.
$$
For further calculation we need the principal part of Laurent series expansion of~$S_0$. In order to find it, consider the function
$$
F(\rho)=\frac{1}{1-e^{-2\rho}}=\frac{1}{\sum_{j=1}^\infty(-1)^{j-1}\frac{(2\rho)^j}{j!}}.
$$
It has a simple pole in $\rho=0$, and its residue equals
$$
\text{Res}(F,0)=\lim_{\rho\to0}\rho\cdot F(\rho)=\frac{1}{2}.
$$
Further,
$$
F(\rho)-\frac{1}{2\rho}=\frac{1}{2}+\frac{\rho}{6}-\frac{\rho^3}{90}+\mathcal O(\rho^4).
$$
Consequently, Laurent series expansion of $S_0(\rho)=\left[F(\rho)\right]^{n-1}$ is in limit $\rho\to0$ equal to
\begin{equation}\label{eq:Laurent_series_S0}
S_0(\rho)=\frac{1}{2^n}\cdot\left[\frac{2}{\rho^{n-1}}+\frac{2\,(n-1)}{\rho^{n-2}}+\frac{(n-1)(n-\tfrac{4}{3})}{\rho^{n-3}}+\frac{\tfrac{1}{3}\,(n-1)^2(n-2)}{\rho^{n-4}}\right]
   +\mathcal O(\rho^{-(n-5)}).
\end{equation}
Moreover, Laurent series of~$S_0$ is convergent for $0<\rho<\infty$ and the relation
$$
S_{m+1}(\rho)=-\frac{1}{2}\,S_m^\prime(\rho)
$$
holds. In order to compute var$_S(f_\rho^m)$ with accuracy of two powers of~$\rho$, we need to develop the numerator and denominator in~\eqref{eq:var_S} with accuracy of four powers of~$\rho$. Thus, since the derivative of the $n$-th term  in~\eqref{eq:Laurent_series_S0} vanishes, the cases $n=4$, $n=3$, and $n=2$ will be treated separately.

Suppose, $n\geq5$. We obtain by induction with respect to~$m$
\begin{equation}\label{eq:Smrho}\begin{split}
S_m(\rho)&=\frac{1}{2^{n+m}}\cdot\left[\frac{2\,(n+m-2)!}{(n-2)!\,\rho^{n+m-1}}+\frac{2\,(n-1)(n+m-3)!}{(n-3)!\,\rho^{n+m-2}}\right.\\
&\left.+\frac{(n-1)(n-\tfrac{4}{3})(n+m-4)!}{(n-4)!\,\rho^{n+m-3}}+\frac{\tfrac{1}{3}\,(n-1)^2(n-2)(n+m-5)!}{(n-5)!\,\rho^{n+m-4}}\right]\\
&+\mathcal{O}(\rho^{-n-m+5})\qquad\text{for }\rho\to0.
\end{split}\end{equation}

Therefore, the variance in space domain~var$_S(f_\rho^m)$ behaves in limit $\rho\to0$ like
\begin{equation}\label{eq:varS_AB}
\text{var}_S(f_\rho^m)=\left(\frac{e^\rho\cdot A}{2\cdot B}\right)^2-1,
\end{equation}
where
\begin{equation}\label{eq:An5}\begin{split}
2&^{n+2m}\cdot A
=\frac{2\,(n+2m-1)!}{(n-1)!\,\rho^{n+2m}}+\frac{2\,(n-1)(n+2m-2)!}{(n-2)!\,\rho^{n+2m-1}}+\frac{(n^2-\tfrac{7}{3}\,n+2)(n+2m-3)!}{(n-3)!\,\rho^{n+2m-2}}\\
&+\frac{\tfrac{1}{3}\,(n-1)(n^2-3n+4)(n+2m-4)!}{(n-4)!\,\rho^{n+2m-3}}+\mathcal{O}(\rho^{-n-2m+4})
\end{split}\end{equation}
(the above formula is derived in Appendix) and
\begin{align*}
B&=\sum_{j=0}^m\binom{m}{j}\left(\frac{S_{m+j+1}}{n-1}+S_{m+j}\right)\\
&=\frac{S_{2m+1}}{n-1}+S_{2m}+m\,\left(\frac{S_{2m}}{n-1}+S_{2m-1}\right)+\frac{m\,(m-1)}{2}\left(\frac{S_{2m-1}}{n-1}+S_{2m-2}\right)\\
&+\frac{m\,(m-1)(m-2)}{6}\left(\frac{S_{2m-2}}{n-1}+\mathcal{O}(\rho^{-n-2m+4})\right)+\mathcal{O}(\rho^{-n-2m+4})\\
&=\frac{S_{2m+1}}{n-1}+\frac{n+m-1}{n-1}\,S_{2m}+\frac{m\,(2n+m-3)}{2\,(n-1)}\,S_{2m-1}\\
&+\frac{m\,(m-1)(3n+m-5)}{6\,(n-1)}\,S_{2m-2}+\mathcal{O}(\rho^{-n-2m+4})
\end{align*}
for $m\geq1$ (in the cases $m=1$ and $m=2$ the rest term vanishes; similarly, some of the series~$S_j$ are multiplied by~$0$). Substituting~\eqref{eq:Smrho} for the series~$S_j$, we obtain
\begin{equation}\label{eq:Bn5}\begin{split}
2&^{n+2m}\cdot B=\frac{(n+2m-1)!}{(n-1)!\,\rho^{n+2m}}+\frac{(n^2-n+2m)(n+2m-2)!}{(n-1)!\,\rho^{n+2m-1}}\\
&+\frac{\left[\tfrac{1}{2}\,n^4-\tfrac{5}{3}\,n^3+(2m+\tfrac{3}{2})\,n^2-\tfrac{1}{3}\,(6m+1)\,n+2m\,(m-1)\right](n+2m-3)!}{(n-1)!\,\rho^{n+2m-2}}\\
&+\left[\tfrac{1}{6}\,n^6-\tfrac{7}{6}\,n^5+(m+\tfrac{17}{6})\,n^4-\tfrac{1}{6}\,(20m+17)\,n^3+(2m^2+m+1)\,n^2\right.\\
&\left.-\tfrac{2}{3}(3m-2)\,mn+\tfrac{4}{3}\,(m^2-3m+2)\,m\right]\cdot\frac{(n+2m-4)!}{(n-1)!\,\rho^{n+2m-3}}+\mathcal{O}(\rho^{-n-2m+4})
\end{split}\end{equation}
(see Appendix for the derivation). In order to compute Taylor series of var$_S(f_\rho^m)$, we multiply both numerator and denominator of the fraction in~\eqref{eq:varS_AB} by~$(2\rho)^{n+2m}$. Taking into account that
$$
e^\rho=1+\rho+\frac{\rho^2}{2}+\frac{\rho^3}{6}+\mathcal{O}(\rho^3),
$$
we obtain
\begin{align*}
e^\rho&\cdot A\cdot(2\rho)^{n+2m}=\frac{2\,(n+2m-1)!}{(n-1)!}+\frac{2\,(n^2-n+2m)(n+2m-2)!}{(n-1)!}\cdot\rho\\
&+\frac{\left[n^4-\tfrac{10}{3}n^3+4\,(m+1)\,n^2-(4m+\tfrac{11}{3})\,n+2\,(2m^2-m+1)\right](n+2m-3)!}{(n-1)!}\cdot\rho^2\\
&+\left[\tfrac{1}{3}\,n^6-\tfrac{7}{3}\,n^5+(2m+\tfrac{20}{3})\,n^4-\tfrac{5}{3}\,(4m+7)\,n^3+(4m^2+6m+13)\,n^2\right.\\
&\left.-2\,(2m^2+\tfrac{8}{3}\,m+3)\,n+\tfrac{4}{3}m\,(2m^2-3m+4)\right]\cdot\frac{(n+2m-4)!}{(n-1)!}\cdot\rho^3+\mathcal O(\rho^4).
\end{align*}
Since
\begin{align*}
&\frac{a+b\rho+c\rho^2+d\rho^3+\mathcal{O}(\rho^4)}{\alpha+\beta\rho+\gamma\rho^2+\delta\rho^3+\mathcal{O}(\rho^4)}
   =\frac{a}{\alpha}+\frac{(\alpha b-a\beta)\,\rho}{\alpha^2}+\frac{(-\alpha b\beta+a\beta^2+\alpha^2c-a\alpha\gamma)\,\rho^2}{\alpha^3}\\
&\quad+\frac{(\alpha^3d-\alpha^2\beta c+\alpha b\beta^2-a\beta^3-\alpha^2b\gamma+2a\alpha\beta\gamma-a\alpha^2\delta)\,\rho^3}{\alpha^4}+\mathcal{O}(\rho^4),
\end{align*}
Taylor series expansion of $\frac{e^\rho\cdot A}{2\cdot B}$ reads
\begin{align*}
\frac{e^\rho\cdot A}{2\cdot B}&=\frac{e^\rho\cdot A\cdot(2\rho)^{n+2m}}{2\cdot B\cdot(2\rho)^{n+2m}}
   =1+\frac{n^2-3n+2m+2}{2\,(n+2m-1)(n+2m-2)}\cdot\rho^2\\
&-\frac{2\,(n-1)^2(n-3)\,m}{(n+2m-1)^2(n+2m-2)(n+2m-3)}\cdot\rho^3+\mathcal{O}(\rho^4).
\end{align*}

Consequently, according to~\eqref{eq:varS_AB},
\begin{align}
\text{var}_S(f_\rho^m)&=\frac{n^2-3n+2m+2}{(n+2m-1)(n+2m-2)}\cdot\rho^2\notag\\
&-\frac{4\,(n-1)^2(n-3)\,m}{(n+2m-1)^2(n+2m-2)(n+2m-3)}\cdot\rho^3+\mathcal{O}(\rho^4).\label{eq:varSn5}
\end{align}

If $n=4$, then
$$
S_m(\rho)=\frac{1}{2^{m+4}}\cdot\left[\frac{(m+2)!}{\rho^{m+3}}+\frac{6\,(m+1)!}{\rho^{m+2}}+\frac{8\,m!}{\rho^{m+1}}\right]
   +\mathcal{O}(1)\qquad\text{for }\rho\to0.
$$
Thus, for $\rho\to0$,
$$
\text{var}_S(f_\rho^m)=\left(\frac{e^\rho\cdot A}{2\cdot B}\right)^2-1
$$
with
\begin{align*}
2&^{2m+4}\cdot A=2^{2m+4}\cdot\left(\frac{2}{3}\,S_{2m+1}+S_{2m}\right)\\
&=\frac{1}{3}\cdot\left[\frac{(2m+3)!}{\rho^{2m+4}}+\frac{6\,(2m+2)!}{\rho^{2m+3}}+\frac{8\,(2m+1)!}{\rho^{2m+2}}\right]\\
&+\frac{(2m+2)!}{\rho^{2m+3}}+\frac{6\,(2m+1)!}{\rho^{2m+2}}+\frac{8\,(2m)!}{\rho^{2m+1}}+\mathcal O(1)\\
&=\frac{\tfrac{1}{3}\,(2m+3)!}{\rho^{2m+4}}+\frac{3\,(2m+2)!}{\rho^{2m+3}}+\frac{\tfrac{26}{3}\,(2m+1)!}{\rho^{2m+2}}+\frac{8\,(2m)!}{\rho^{2m+1}}
   +\mathcal O(1)
\end{align*}
and
\begin{align*}
2^{2m+4}\cdot B&=2^{2m+4}\cdot\left[\frac{S_{2m+1}}{3}+\frac{m+3}{3}\,S_{2m}+\frac{m\,(m+5)}{6}\,S_{2m-1}\right.\\
&\left.+\frac{m\,(m-1)(m+8)}{18}\,S_{2m-2}\right]+\mathcal{O}(\rho^{-2m})\\
&=\frac{1}{6}\cdot\left[\frac{(2m+3)!}{\rho^{2m+4}}+\frac{6\,(2m+2)!}{\rho^{2m+3}}+\frac{8\,(2m+1)!}{\rho^{2m+2}}\right]\\
&+\frac{m+3}{3}\cdot\left[\frac{(2m+2)!}{\rho^{2m+3}}+\frac{6\,(2m+1)!}{\rho^{2m+2}}+\frac{8\,(2m)!}{\rho^{2m+1}}\right]\\
&+\frac{m\,(m+5)}{3}\cdot\left[\frac{(2m+1)!}{\rho^{2m+2}}+\frac{6\,(2m)!}{\rho^{2m+1}}\right]+\frac{2m\,(m-1)(m+8)}{9}\cdot\frac{(2m)!}{\rho^{2m+1}}
   +\mathcal{O}(\rho^{-2m})\\
&=\frac{\tfrac{1}{6}\cdot(2m+3)!}{\rho^{2m+4}}+\frac{\tfrac{m+6}{3}\cdot(2m+2)!}{\rho^{2m+3}}+\frac{\tfrac{m^2+11m+22}{3}\cdot(2m+1)!}{\rho^{2m+2}}\\
&+\frac{\tfrac{2\,(m^3+16m^2+49m+36)}{9}\cdot(2m)!}{\rho^{2m+1}}+\mathcal{O}(\rho^{-2m}).
\end{align*}
Therefore,
\begin{align*}
e^\rho&\cdot A\cdot(2\rho)^{2m+4}=\frac{2\,(m+1)(2m+3)(2m+1)!}{3}+\frac{4\,(m+1)(m+6)(2m+1)!}{3}\cdot\rho\\
&+\frac{(2m^2+23m+47)(2m+1)!}{3}\cdot\rho^2+\frac{2}{9}\,(2m^3+33m^2+124m+90)(2m)!\cdot\rho^3+\mathcal O(\rho^4).
\end{align*}
Further,
$$
\frac{e^\rho\cdot A}{2\cdot B}=1+\frac{(m+3)}{2\,(2m^2+5m+3)}\cdot\rho^2-\frac{m\,(4m^2+2m+21)}{3\,(2m+3)^2(2m^2+3m+1)}\cdot\rho^3+\mathcal O(\rho^4)
$$
and
$$
\text{var}_S(f_\rho^m)=\frac{(m+3)}{2m^2+5m+3}\cdot\rho^2-\frac{2m\,(4m^2+2m+21)}{3\,(2m+3)^2(2m^2+3m+1)}\cdot\rho^3+\mathcal O(\rho^4).
$$

In the case $n=3$,
$$
S_m(\rho)=\frac{1}{2^{m+2}}\cdot\left[\frac{(m+1)!}{\rho^{m+2}}+\frac{2\,m!}{\rho^{m+1}}\right]+\mathcal{O}(1)\qquad\text{for }\rho\to0.
$$
Consequently, for $\rho\to0$,
$$
\text{var}_S(f_\rho^m)=\left(\frac{e^\rho\cdot A}{2\cdot B}\right)^2-1,
$$
where
\begin{align*}
2&^{2m+3}\cdot A=2^{2m+3}\cdot\left(S_{2m+1}+S_{2m}\right)\\
&=\frac{(2m+2)!}{\rho^{2m+3}}+\frac{2\,(2m+1)!}{\rho^{2m+2}}+\frac{2\,(2m+1)!}{\rho^{2m+2}}+\frac{4\,(2m)!}{\rho^{2m+1}}+\mathcal O(1)\\
&=\frac{(2m+2)!}{\rho^{2m+3}}+\frac{4\,(2m+1)!}{\rho^{2m+2}}+\frac{4\,(2m)!}{\rho^{2m+1}}+\mathcal O(1)
\end{align*}
and
\begin{align*}
2&^{2m+3}\cdot B=\frac{S_{2m+1}}{2}+\frac{m+2}{2}\,S_{2m}+\frac{m\,(m+3)}{4}\,S_{2m-1}\\
&+\frac{m\,(m-1)(m+4)}{12}\,S_{2m-2}+\mathcal{O}(\rho^{-2m+1})\\
&=\frac{1}{2}\cdot\left[\frac{(2m+2)!}{\rho^{2m+3}}+\frac{2\,(2m+1)!}{\rho^{2m+2}}\right]+(m+2)\left[\frac{(2m+1)!}{\rho^{2m+2}}+\frac{2\,(2m)!}{\rho^{2m+1}}\right]\\
&+m(m+3)\left[\frac{(2m)!}{\rho^{2m+1}}+\frac{2\,(2m-1)!}{\rho^{2m}}\right]+\frac{2m\,(m-1)(m+4)}{3}\cdot\frac{(2m-1)!}{\rho^{2m}}+\mathcal O(\rho^{-2m+1})\\
&=\frac{\tfrac{1}{2}\,(2m+2)!}{\rho^{2m+3}}+\frac{(m+3)(2m+1)!}{\rho^{2m+2}}+\frac{(m+1)(m+4)\,(2m)!}{\rho^{2m+1}}\\
&+\frac{\tfrac{2}{3}\,m\,(m+1)(m+5)(2m-1)!}{\rho^{2m}}+\mathcal O(\rho^{-2m+1})
\end{align*}
Therefore,
\begin{align*}
e^\rho&\cdot A\cdot(2\rho)^{2m+3}=(2m+2)!+2\,(m+3)(2m+1)!\cdot\rho\\
&+(m+1)(2m+9)(2m)!\cdot\rho^2+\frac{(2m^2+15m+19)(2m)!}{3}\cdot\rho^3+\mathcal O(\rho^4)
\end{align*}
and
$$
\frac{e^\rho\cdot A}{2\cdot B}=1+\frac{1}{2\,(2m+1)}\cdot\rho^2-\frac{(m-1)(m+1)(m+5)}{3\,(2m+1)}\cdot\rho^3+\mathcal O(\rho^4).
$$
Hence,
$$
\text{var}_S(f_\rho^m)=\frac{1}{2m+1}\cdot\rho^2-\frac{2\,(m-1)(m+1)(m+5)}{3\,(2m+1)}\cdot\rho^3+\mathcal O(\rho^4).
$$

Finally, if $n=2$, then
$$
S_m(\rho)=\frac{1}{2^{m+1}}\cdot\frac{m!}{\rho^{m+1}}+\mathcal{O}(1)\qquad\text{for }\rho\to0.
$$
Hence, for $\rho\to0$,
\begin{align*}
\text{var}&_S(f_\rho^m)=\left[\frac{e^\rho\cdot(2\,S_{2m+1}+S_{2m})}{2\cdot\sum_{j=1}^{m+1}\binom{m+1}{j}S_{m+j}}\right]^2-1\\
&=\left[\frac{\left(1+\frac{2\rho}{2}+\frac{(2\rho)^2}{8}+\frac{(2\rho)^3}{48}+\mathcal O(\rho^4)\right)
   \overbrace{\left(\tfrac{2\,(2m+1)!}{(2\rho)^{2m+2}}+\tfrac{(2m)!}{(2\rho)^{2m+1}}+\mathcal O(1)\right)}^{A}}
   {2\cdot\left(\frac{(2m+1)!}{(2\rho)^{2m+2}}+\frac{(m+1)(2m)!}{(2\rho)^{2m+1}}+\frac{(m+1)\,m\,(2m-1)!}{2\,(2\rho)^{2m}}
   +\frac{(m+1)\,m\,(m-1)(2m-2)!}{6\,(2\rho)^{2m-1}}\right)+\mathcal O(\rho^{-2m+2})}\right]^2-1\\
&=\left[\frac{\frac{2\,(2m+1)!}{(2\rho)^{2m+2}}+\frac{2\,(m+1)(2m)!}{(2\rho)^{2m+1}}+\frac{(2m+3)(2m)!}{4\,(2\rho)^{2m}}
   +\frac{(m+2)(2m)!}{12\,(2\rho)^{2m-1}}+\mathcal O(\rho^{-2m+2})}
   {2\cdot\left(\frac{(2m+1)!}{(2\rho)^{2m+2}}+\frac{(m+1)(2m)!}{(2\rho)^{2m+1}}+\frac{(m+1)(2m)!}{4\,(2\rho)^{2m}}
   +\frac{(m+1)\,m\,(m-1)(2m-2)!}{6\,(2\rho)^{2m-1}}\right)+\mathcal O(\rho^{-2m+2})}\right]^2-1\\
&=\left[1+\frac{\rho^2}{2\,(2m+1)}+\frac{\rho^3}{(2m-1)(2m+1)^2}+\mathcal O(\rho^4)\right]^2-1\\
&=\frac{\rho^2}{2m+1}+\frac{2\rho^3}{(2m-1)(2m+1)^2}+\mathcal O(\rho^4).
\end{align*}
Note that this formula matches with~\eqref{eq:varSn5} with $n=2$ substituted.

On the other hand,
$$
\text{var}_M(f_\rho^m)=\frac{C}{A}
$$
for
\begin{equation}\label{eq:varMn3}
2^{n+2m+2}\cdot C=\frac{2\,(n+2m+1)!}{(n-1)!\,\rho^{n+2m+2}}+\frac{2\,(n+1)\,(n+2m)!}{(n-2)!\,\rho^{n+2m+1}}+\mathcal{O}(\rho^{-n-2m})
\end{equation}
if $n\geq3$ (see Appendix), respectively
\begin{align*}
2^{2m+4}\cdot C&=2^{2m+4}\cdot\left(2\,S_{2m+3}+3\,S_{2m+2}+S_{2m+1}\right)\\
&=\frac{2\,(2m+3)!}{\rho^{2m+4}}+\frac{6\,(2m+2)!}{\rho^{2m+3}}+\mathcal O(\rho^{-2m-2})
\end{align*}
if $n=2$, i.e., formula~\eqref{eq:varMn3} can be used also in this case. Similarly, the first two elements in the series representing~$A$ match with those given for $n\geq5$ in formula~\eqref{eq:An5}.
 Thus,
$$
\text{var}_M(f_\rho^m)=\frac{(n+2m)(n+2m+1)}{4\rho^2}+\frac{(n-1)\,m\,(n+2m)}{(n+2m-1)\,\rho}+\mathcal O(1)
$$
for $n\geq2$ and $\rho\to0$. The values of the uncertainty product are obtained by multiplication of square roots of var$_M(f_\rho^m)$ and var$_S(f_\rho^m)$.\hfill$\Box$\\

Since only two terms of the series expansion of the uncertainty product of~$g_\rho^m$ are derived, it is impossible to find the minimum of~$U(g_\rho^m)$ with respect to~$\rho$. However, in analogy to the problem studied in~\cite{FP02}, consider the limit value of the uncertainty product for $\rho\to0$,
$$
\lim_{\rho\to0}U(g_\rho^m)=\begin{cases}\frac{1}{2}\,\sqrt{\frac{(n+2m)(n+2m+1)[n^2-3n+2\,(m+1)]}{(n+2m-1)(n+2m-2)}},&n\geq5,\,n=2,\\
\sqrt{\frac{(m+3)(m+2)(2m+5)}{2(m+1)(2m+3)}},&n=4,\\
\sqrt{\frac{(m+2)(2m+3)}{2(2m+1)}},&n=3.\end{cases}
$$

It is interesting to observe how it behaves. It can be easily verified that it is increasing in~$m$ for $n=2,\,3,\,4$. In the case $n\geq5$ consider the function
$$
F:\,m\mapsto\frac{(n+2m)(n+2m+1)[n^2-3n+2\,(m+1)]}{(n+2m-1)(n+2m-2)},
$$
defined for $m\in\mathbb R$. Its derivative
$$
F^\prime:\,m\mapsto\tfrac{2[16m^4\!+\!16m^3(2n\!-\!3)\!+\!4m^2(2n^2\!-\!2n\!-\!5)\!-\!4m(2n^3\!-\!9n^2\!+\!13n\!-\!6)\!-\!(n\!-\!2)^2 (3n^2\!-\!2n\!-\!1)]}
    {(n+2m-1)^2\,(n+2m-2)^2}
$$
has four roots. Since
$$
F^\prime(m)
\begin{cases}>0&\text{for }m=-\tfrac{3n}{2},\\<0&\text{for }m=-\tfrac{n}{2}-1,\\>0&\text{for }m=-\tfrac{n}{2},\\<0&\text{for }m=-\frac{n}{2}+\tfrac{3}{4},\end{cases}
$$
three of the roots are negative. Further,
$$
F^\prime(m)
\begin{cases}<0&\text{for }m=\tfrac{n}{2}-1,\\>0&\text{for }m=\tfrac{n}{2}-\tfrac{1}{2},\end{cases}
$$
i.e., function~$F$ has its local minimum in $m\in(\tfrac{n}{2}-1,\tfrac{n}{2}-\tfrac{1}{2})$. Since
$$
\lim_{\rho\to0}U(g_\rho^m)=\frac{1}{2}\sqrt{F(m)}
$$
and
$$
F(\tfrac{n}{2}-1)=F(\tfrac{n}{2}-\tfrac{1}{2})=\frac{n(n-1)(2n-1)}{2n-3},
$$
minimum value of the limit of the uncertainty product is equal to
\begin{equation}\label{eq:min_Ug}
\min_{m\in\mathbb N}\lim_{\rho\to0}U(g_\rho^m)=\frac{1}{2}\sqrt{\frac{n(n-1)(2n-1)}{2n-3}},
\end{equation}
and it is obtained for $m=\left[\frac{n-1}{2}\right]$. If $n\to\infty$, the value of the expression~\eqref{eq:min_Ug} behaves like~$\frac{n}{2}$, i.e. it tends to the minimum value according to~\eqref{eq:UP}. It is well-known that the minimum uncertainty product is a property of Gauss functions, see~\cite{FM99,FP02} for the case of functions over~$\mathcal S^2$. Thus, it is an interesting and surprising statement that another family of functions has a similar feature.

\section{Appendix}

Formula~\eqref{eq:An5} is obtained in the following way:\nopagebreak\vspace{-2em}

\begin{align*}
2&^{n+2m}\cdot A=2^{n+2m}\cdot\left(\tfrac{2}{n-1}\,S_{2m+1}+S_{2m}\right)\\[-0.2em]
&=\tfrac{1}{n-1}\cdot\left[\tfrac{2\,(n+2m-1)!}{(n-2)!\,\rho^{n+2m}}+\tfrac{2\,(n-1)(n+2m-2)!}{(n-3)!\,\rho^{n+2m-1}}\right.\\[-0.2em]
&\left.+\tfrac{(n-1)(n-\tfrac{4}{3})(n+2m-3)!}{(n-4)!\,\rho^{n+2m-2}}+\tfrac{(n-1)^2(n-2)(n+2m-4)!}{3\,(n-5)!\,\rho^{n+2m-3}}\right]\\[-0.2em]
&+\tfrac{2\,(n+2m-2)!}{(n-2)!\,\rho^{n+2m-1}}+\tfrac{2\,(n-1)(n+2m-3)!}{(n-3)!\,\rho^{n+2m-2}}
   +\tfrac{(n-1)(n-\tfrac{4}{3})(n+2m-4)!}{(n-4)!\,\rho^{n+2m-3}}+\mathcal{O}(\rho^{-n-2m+4})\\[0.8em]
&=\tfrac{2\,(n+2m-1)!}{(n-1)!\,\rho^{n+2m}}+\tfrac{2\,(n+2m-2)!}{(n-3)!\,\rho^{n+2m-1}}
   +\tfrac{(n-\tfrac{4}{3})(n+2m-3)!}{(n-4)!\,\rho^{n+2m-2}}+\tfrac{(n-1)(n-2)(n+2m-4)!}{3\,(n-5)!\,\rho^{n+2m-3}}\\[-0.2em]
&+\tfrac{2\,(n+2m-2)!}{(n-2)!\,\rho^{n+2m-1}}+\tfrac{2\,(n-1)(n+2m-3)!}{(n-3)!\,\rho^{n+2m-2}}
   +\tfrac{(n-1)(n-\tfrac{4}{3})(n+2m-4)!}{(n-4)!\,\rho^{n+2m-3}}+\mathcal{O}(\rho^{-n-2m+4})\\[0.8em]
&=\tfrac{2\,(n+2m-1)!}{(n-1)!\,\rho^{n+2m}}+\tfrac{2\,(n+2m-2)!}{(n-3)!\,\rho^{n+2m-1}}+\tfrac{2\,(n+2m-2)!}{(n-2)!\,\rho^{n+2m-1}}\\[-0.2em]
&+\tfrac{(n-\tfrac{4}{3})(n+2m-3)!}{(n-4)!\,\rho^{n+2m-2}}+\tfrac{2\,(n-1)(n+2m-3)!}{(n-3)!\,\rho^{n+2m-2}}\\[-0.2em]
&+\tfrac{(n-1)(n-2)(n+2m-4)!}{3\,(n-5)!\,\rho^{n+2m-3}}+\tfrac{(n-1)(n-\tfrac{4}{3})(n+2m-4)!}{(n-4)!\,\rho^{n+2m-3}}+\mathcal{O}(\rho^{-n-2m+4})\\[0.8em]
&=\tfrac{2\,(n+2m-1)!}{(n-1)!\,\rho^{n+2m}}+\tfrac{2\,(n-1)(n+2m-2)!}{(n-2)!\,\rho^{n+2m-1}}
   +\tfrac{\left[(n-\tfrac{4}{3})(n-3)+2\,(n-1)\right](n+2m-3)!}{(n-3)!\,\rho^{n+2m-2}}\\[-0.2em]
&+\tfrac{(n-1)\left[(n-2)(n-4)+3\,(n-\tfrac{4}{3})\right](n+2m-4)!}{3\,(n-4)!\,\rho^{n+2m-3}}+\mathcal{O}(\rho^{-n-2m+4})\\[0.8em]
&=\tfrac{2\,(n+2m-1)!}{(n-1)!\,\rho^{n+2m}}+\tfrac{2\,(n-1)(n+2m-2)!}{(n-2)!\,\rho^{n+2m-1}}
   +\tfrac{(n^2-\tfrac{7}{3}\,n+2)(n+2m-3)!}{(n-3)!\,\rho^{n+2m-2}}\\[-0.2em]
&+\tfrac{\tfrac{1}{3}\,(n-1)(n^2-3n+4)(n+2m-4)!}{(n-4)!\,\rho^{n+2m-3}}+\mathcal{O}(\rho^{-n-2m+4}).
\end{align*}

Formula~\eqref{eq:Bn5} is obtained in the following way:\nopagebreak\vspace{-2em}

\begin{align*}
2&^{n+2m}\cdot B=\tfrac{1}{2\,(n-1)}\cdot\left[\tfrac{2\,(n+2m-1)!}{(n-2)!\,\rho^{n+2m}}+\tfrac{2\,(n-1)(n+2m-2)!}{(n-3)!\,\rho^{n+2m-1}}\right.\\[-0.2em]
&\left.+\tfrac{(n-1)(n-\tfrac{4}{3})(n+2m-3)!}{(n-4)!\,\rho^{n+2m-2}}+\tfrac{(n-1)^2(n-2)(n+2m-4)!}{3\,(n-5)!\,\rho^{n+2m-3}}\right]\\[-0.2em]
&+\tfrac{n+m-1}{n-1}\cdot\left[\tfrac{2\,(n+2m-2)!}{(n-2)!\,\rho^{n+2m-1}}+\tfrac{2\,(n-1)(n+2m-3)!}{(n-3)!\,\rho^{n+2m-2}}
   +\tfrac{(n-1)(n-\tfrac{4}{3})(n+2m-4)!}{(n-4)!\,\rho^{n+2m-3}}\right]\\[-0.2em]
&+\tfrac{m\,(2n+m-3)}{(n-1)}\cdot\left[\tfrac{2\,(n+2m-3)!}{(n-2)!\,\rho^{n+2m-2}}+\tfrac{2\,(n-1)(n+2m-4)!}{(n-3)!\,\rho^{n+2m-3}}\right]\\[-0.2em]
&+\tfrac{2m\,(m-1)(3n+m-5)}{3\,(n-1)}\cdot\tfrac{2\,(n+2m-4)!}{(n-2)!\,\rho^{n+2m-3}}+\mathcal{O}(\rho^{-n-2m+4})\\[0.8em]
&=\tfrac{(n+2m-1)!}{(n-1)!\,\rho^{n+2m}}+\tfrac{(n+2m-2)!}{(n-3)!\,\rho^{n+2m-1}}+\tfrac{(n-\tfrac{4}{3})(n+2m-3)!}{2\,(n-4)!\,\rho^{n+2m-2}}
   +\tfrac{(n-1)(n-2)(n+2m-4)!}{6\,(n-5)!\,\rho^{n+2m-3}}\\[-0.2em]
&+\tfrac{2\,(n+m-1)(n+2m-2)!}{(n-1)!\,\rho^{n+2m-1}}+\tfrac{2\,(n+m-1)(n+2m-3)!}{(n-3)!\,\rho^{n+2m-2}}
   +\tfrac{(n+m-1)(n-\tfrac{4}{3})(n+2m-4)!}{(n-4)!\,\rho^{n+2m-3}}\\[-0.2em]
&+\tfrac{2m\,(2n+m-3)(n+2m-3)!}{(n-1)!\,\rho^{n+2m-2}}+\tfrac{2m\,(2n+m-3)(n+2m-4)!}{(n-3)!\,\rho^{n+2m-3}}\\[-0.2em]
&+\tfrac{4m\,(m-1)(3n+m-5)(n+2m-4)!}{3\,(n-1)!\,\rho^{n+2m-3}}+\mathcal{O}(\rho^{-n-2m+4})\\[0.8em]
&=\tfrac{(n+2m-1)!}{(n-1)!\,\rho^{n+2m}}+\tfrac{(n+2m-2)!}{(n-3)!\,\rho^{n+2m-1}}+\tfrac{2\,(n+m-1)(n+2m-2)!}{(n-1)!\,\rho^{n+2m-1}}\\[-0.2em]
&+\tfrac{(n-\tfrac{4}{3})(n+2m-3)!}{2\,(n-4)!\,\rho^{n+2m-2}}+\tfrac{2\,(n+m-1)(n+2m-3)!}{(n-3)!\,\rho^{n+2m-2}}
   +\tfrac{2m\,(2n+m-3)(n+2m-3)!}{(n-1)!\,\rho^{n+2m-2}}\\[-0.2em]
&+\tfrac{(n-1)(n-2)(n+2m-4)!}{6\,(n-5)!\,\rho^{n+2m-3}}+\tfrac{(n+m-1)(n-\tfrac{4}{3})(n+2m-4)!}{(n-4)!\,\rho^{n+2m-3}}
   +\tfrac{2m\,(2n+m-3)(n+2m-4)!}{(n-3)!\,\rho^{n+2m-3}}\\[-0.2em]
&+\tfrac{4m\,(m-1)(3n+m-5)(n+2m-4)!}{3\,(n-1)!\,\rho^{n+2m-3}}+\mathcal{O}(\rho^{-n-2m+4})\\[0.2em]
&=\tfrac{(n+2m-1)!}{(n-1)!\,\rho^{n+2m}}+\tfrac{\left[(n-1)(n-2)+2\,(n+m-1)\right](n+2m-2)!}{(n-1)!\,\rho^{n+2m-1}}\\[-0.2em]
&+\tfrac{\left[\tfrac{1}{2}(n-\tfrac{4}{3})(n-1)(n-2)(n-3)+2\,(n+m-1)(n-1)(n-2)+2m\,(2n+m-3)\right](n+2m-3)!}{(n-1)!\,\rho^{n+2m-2}}\\[-0.2em]
&+\tfrac{(n+2m-4)!}{(n-1)!\,\rho^{n+2m-3}}\cdot\left[\tfrac{1}{6}\,(n-1)^2(n-2)^2(n-3)(n-4)\right.\\[-0.2em]
&+(n+m-1)(n-\tfrac{4}{3})(n-1)(n-2)(n-3)\\[-0.2em]
&\left.2m\,(2n+m-3)(n-1)(n-2)+\tfrac{4}{3}\,m\,(m-1)(3n+m-5)\right]+\mathcal{O}(\rho^{-n-2m+4})\\[0.8em]
&=\tfrac{(n+2m-1)!}{(n-1)!\,\rho^{n+2m}}+\tfrac{(n^2-n+2m)(n+2m-2)!}{(n-1)!\,\rho^{n+2m-1}}\\[-0.2em]
&+\tfrac{\left[\tfrac{1}{2}\,n^4-\tfrac{5}{3}\,n^3+(2m+\tfrac{3}{2})\,n^2-\tfrac{1}{3}\,(6m+1)\,n+2m\,(m-1)\right](n+2m-3)!}{(n-1)!\,\rho^{n+2m-2}}\\[-0.2em]
&+\left[\tfrac{1}{6}\,n^6-\tfrac{7}{6}\,n^5+(m+\tfrac{17}{6})\,n^4-\tfrac{1}{6}\,(20m+17)\,n^3+(2m^2+m+1)\,n^2\right.\\[-0.2em]
&\left.-\tfrac{2}{3}(3m-2)\,mn+\tfrac{4}{3}\,(m^2-3m+2)\,m\right]\cdot\tfrac{(n+2m-4)!}{(n-1)!\,\rho^{n+2m-3}}+\mathcal{O}(\rho^{-n-2m+4}).
\end{align*}

Formula~\eqref{eq:varMn3} is obtained in the following way:\nopagebreak\vspace{-2em}

\begin{align*}
2&^{n+2m}\cdot C=\tfrac{1}{4\,(n-1)}\cdot\left[\tfrac{2\,(n+2m+1)!}{(n-2)!\,\rho^{n+2m+2}}+\tfrac{2\,(n-1)(n+2m)!}{(n-3)!\,\rho^{n+2m+1}}\right.\\[-0.2em]
&\left.+\tfrac{(n-1)(n-\tfrac{4}{3})(n+2m-1)!}{(n-4)!\,\rho^{n+2m}}+\tfrac{(n-1)^2(n-2)(n+2m-2)!}{3\,(n-5)!\,\rho^{n+2m-1}}\right]\\[-0.2em]
&+\tfrac{3}{4}\cdot\left[\tfrac{2\,(n+2m)!}{(n-2)!\,\rho^{n+2m+1}}+\tfrac{2\,(n-1)(n+2m-1)!}{(n-3)!\,\rho^{n+2m}}
   +\tfrac{(n-1)(n-\tfrac{4}{3})(n+2m-2)!}{(n-4)!\,\rho^{n+2m-1}}\right]\\[-0.2em]
&+\tfrac{n-1}{2}\cdot\left[\tfrac{2\,(n+2m-1)!}{(n-2)!\,\rho^{n+2m}}+\tfrac{2\,(n-1)(n+2m-2)!}{(n-3)!\,\rho^{n+2m-1}}\right]+\mathcal{O}(\rho^{-n-2m+1})\\[0.8em]
&=\tfrac{(n+2m+1)!}{2\,(n-1)!\,\rho^{n+2m+2}}+\tfrac{(n+2m)!}{2\,(n-3)!\,\rho^{n+2m+1}}+\tfrac{(n-\tfrac{4}{3})(n+2m-1)!}{4\,(n-4)!\,\rho^{n+2m}}\\[-0.2em]
&+\tfrac{(n-1)(n-2)(n+2m-2)!}{12\,(n-5)!\,\rho^{n+2m-1}}\\[-0.2em]
&+\tfrac{3\,(n+2m)!}{2\,(n-2)!\,\rho^{n+2m+1}}+\tfrac{3\,(n-1)(n+2m-1)!}{2\,(n-3)!\,\rho^{n+2m}}+\tfrac{3\,(n-1)(n-\tfrac{4}{3})(n+2m-2)!}{4\,(n-4)!\,\rho^{n+2m-1}}\\[-0.2em]
&+\tfrac{(n-1)(n+2m-1)!}{(n-2)!\,\rho^{n+2m}}+\tfrac{(n-1)^2(n+2m-2)!}{(n-3)!\,\rho^{n+2m-1}}+\mathcal{O}(\rho^{-n-2m+1})\\[0.8em]
&=\tfrac{(n+2m+1)!}{2\,(n-1)!\,\rho^{n+2m+2}}+\tfrac{(n+2m)!}{2\,(n-3)!\,\rho^{n+2m+1}}+\tfrac{3\,(n+2m)!}{2\,(n-2)!\,\rho^{n+2m+1}}\\[-0.2em]
&+\tfrac{(n-\tfrac{4}{3})(n+2m-1)!}{4\,(n-4)!\,\rho^{n+2m}}+\tfrac{3\,(n-1)(n+2m-1)!}{2\,(n-3)!\,\rho^{n+2m}}+\tfrac{(n-1)(n+2m-1)!}{(n-2)!\,\rho^{n+2m}}\\[-0.2em]
&+\tfrac{(n-1)(n-2)(n+2m-2)!}{12\,(n-5)!\,\rho^{n+2m-1}}+\tfrac{3\,(n-1)(n-\tfrac{4}{3})(n+2m-2)!}{4\,(n-4)!\,\rho^{n+2m-1}}\\[-0.2em]
&+\tfrac{(n-1)^2(n+2m-2)!}{(n-3)!\,\rho^{n+2m-1}}+\mathcal{O}(\rho^{-n-2m+1})\\[0.8em]
&=\tfrac{(n+2m+1)!}{2\,(n-1)!\,\rho^{n+2m+2}}+\tfrac{(n-2+3)\,(n+2m)!}{2\,(n-2)!\,\rho^{n+2m+1}}\\[-0.2em]
&+\tfrac{\left[(n-\tfrac{4}{3})(n-2)(n-3)+6\,(n-1)(n-2)+4\,(n-1)\right](n+2m-1)!}{4\,(n-2)!\,\rho^{n+2m}}\\[-0.2em]
&+\left[(n-2)(n-3)(n-4)+9(n-\tfrac{4}{3})(n-3)+12\,(n-1)\right]\cdot\tfrac{(n-1)(n+2m-2)!}{12\,(n-3)!\,\rho^{n+2m-1}}\\[-0.2em]
&+\mathcal{O}(\rho^{-n-2m+1})\\[0.8em]
&=\tfrac{\tfrac{1}{2}\,(n+2m+1)!}{(n-1)!\,\rho^{n+2m+2}}+\tfrac{\tfrac{1}{2}\,(n+1)\,(n+2m)!}{(n-2)!\,\rho^{n+2m+1}}
   +\tfrac{\tfrac{1}{12}\,n\,(3n^2-n-4)(n+2m-1)!}{(n-2)!\,\rho^{n+2m}}\\[-0.2em]
&+\tfrac{\tfrac{1}{12}\,(n+1)\,n\,(n-1)^2(n+2m-2)!}{(n-3)!\,\rho^{n+2m-1}}+\mathcal{O}(\rho^{-n-2m+1}).
\end{align*}

\end{document}